\documentclass[a4paper,twoside,12pt]{amsart}
\RequirePackage{amsmath,amssymb,xypic}
\theoremstyle{plain}
\newtheorem{theorem}{Theorem}[section]
\newtheorem{corollary}[theorem]{Corollary}
\newtheorem{proposition}[theorem]{Proposition}
\newtheorem{lemma}[theorem]{Lemma}

\theoremstyle{definition}
\newtheorem{definition}[theorem]{Definition}
\newtheorem{example}[theorem]{Example}
\newtheorem{remark}[theorem]{Remark}

\numberwithin{equation}{section}

\newcommand{\Union}{\bigcup\limits}

\newcommand{\C}{\mathbb{C}}
\newcommand{\N}{\mathbb{N}}

\newcommand{\Z}{\mathbb{Z}}

\newcommand{\op}{\mathrm{op}}

\DeclareMathOperator{\id}{id}

\newcommand{\ilim}[1][]{\mathop{\varinjlim}\limits_{#1}}

\renewcommand{\to}[1][]{\xrightarrow[]{#1}}

\newcommand{\isoto}[1][]{\xrightarrow[#1]{\sim}}

\newcommand{\Endo}[1][]{\mathrm{End}_{\raise1.5ex\hbox to.1em{}#1}}

\newcommand{\Hom}[1][]{\mathrm{Hom}_{\raise1.5ex\hbox to.1em{}#1}}

\newcommand{\RHom}[1][]{\mathrm{RHom}_{\raise1.5ex\hbox to.1em{}#1}}

\newcommand{\Ext}[2][]{\mathrm{Ext}_{\raise1.5ex\hbox to.1em{}#1}^{#2}}

\newcommand{\THom}[1][]{\mathrm{THom}_{\raise1.5ex\hbox to.1em{}#1}}

\newcommand{\Tens}[1][]{\mathbin{\otimes_{\raise1.5ex\hbox to-.1em{}#1}}}

\newcommand{\LTens}[1][]{\mathbin{\otimes_{\raise1.5ex\hbox to-.1em{}#1}^{L}}}

\newcommand{\Tor}[2][]{\mathrm{Tor}^{\raise1.5ex\hbox to.1em{}#1}_{#2}}

\def\sha{\mathcal{A}}

\def\she{\mathcal{E}}

\def\shg{\mathcal{G}}

\def\shl{\mathcal{L}}

\def\sho{\mathcal{O}}
\def\shp{\mathcal{P}}
\def\shq{\mathcal{Q}}
\def\shr{\mathcal{R}}

\def\shu{\mathcal{U}}

\def\shw{\mathcal{W}}

\newcommand{\shendo}[1][]{{\mathcal{E}nd}_{\raise1.5ex\hbox to.1em{}#1}}

\newcommand{\shaut}[1][]{{\mathcal{A}ut}_{\raise1.5ex\hbox to.1em{}#1}}

\renewcommand{\hom}[1][]{{\mathcal{H}om}_{\raise1.5ex\hbox to.1em{}#1}}

\newcommand{\rhom}[1][]{{R\mathcal{H}om}_{\raise1.5ex\hbox to.1em{}#1}}

\newcommand{\ext}[2][]{{\mathcal{E}xt}_{\raise1.5ex\hbox to.1em{}#1}^{#2}}

\newcommand{\thom}[1][]{{T\mathcal{H}om}_{\raise1.5ex\hbox to.1em{}#1}}

\newcommand{\tens}[1][]{\mathbin{\otimes_{\raise1.5ex\hbox to-.1em{}#1}}}

\newcommand{\ltens}[1][]{\mathbin{\otimes_{\raise1.5ex\hbox to-.1em{}#1}^{L}}}

\newcommand{\tor}[2][]{{\mathcal{T}or}^{\raise1.5ex\hbox to.1em{}#1}_{#2}}

\newcommand{\oim}[1]{{#1}_*}

\newcommand{\opb}[1]{#1^{-1}}

\newcommand{\GHom}[1][]{\mathrm{GHom}_{\raise1.5ex\hbox to.1em{}#1}}

\newcommand{\GExt}[2][]{\mathrm{GExt}_{\raise1.5ex\hbox to.1em{}#1}^{#2}}

\newcommand{\FHom}[1][]{\mathrm{FHom}_{\raise1.5ex\hbox to.1em{}#1}}

\newcommand{\ghom}[1][]{{\mathcal{GH}om}_{\raise1.5ex\hbox to.1em{}#1}}

\newcommand{\gext}[2][]{{\mathcal{GE}xt}_{\raise1.5ex\hbox to.1em{}#1}^{#2}}

\newcommand{\fhom}[1][]{{\mathcal{FH}om}_{\raise1.5ex\hbox to.1em{}#1}}

\newcommand{\gr}{\mathop{\mathcal{G}r}\nolimits}
\newcommand{\Gr}{\mathop{\mathrm{Gr}}\nolimits}

\newcommand{\tenstop}[1][]{\mathbin{\hat{\otimes}_{\raise1.5ex\hbox to-.1em{}#1}}}

\newcommand{\homtop}[1][]{\mathcal{L}_{\raise1.5ex\hbox to.1em{}#1}}

\newcommand{\Homtop}[1][]{\mathrm{L}_{\raise1.5ex\hbox to.1em{}#1}}

\renewcommand{\O}{\mathcal{O}}

\def\absdoim#1{\underline{#1}_*}
\def\reldoim[#1]#2{\underline{#2}_{|{#1}*}}
\def\doim{\@ifnextchar [{\reldoim}{\absdoim}}

\def\absdeim#1{\underline{#1}_*}
\def\reldeim[#1]#2{\underline{#2}_{|{#1}*}}
\def\deim{\@ifnextchar [{\reldeim}{\absdeim}}

\def\absdopb#1{\underline{#1}^{-1}}
\def\reldopb[#1]#2{\underline{#2}_{|{#1}}^{-1}}
\def\dopb{\@ifnextchar [{\reldopb}{\absdopb}}

\def\absboim#1{\underline{\underline{#1}}_*}
\def\relboim[#1]#2{\underline{\underline{#2}}_{|{#1}*}}
\def\boim{\@ifnextchar [{\relboim}{\absboim}}

\def\absbeim#1{\underline{\underline{#1}}_*}
\def\relbeim[#1]#2{\underline{\underline{#2}}_{|{#1}*}}
\def\beim{\@ifnextchar [{\relbeim}{\absbeim}}

\def\absbopb#1{\underline{\underline{#1}}^*}
\def\relbopb[#1]#2{\underline{\underline{#2}}_{|{#1}}^*}
\def\bopb{\@ifnextchar [{\relbopb}{\absbopb}}

\newcommand{\ad}{\operatorname{ad}}
\newcommand{\WKB}{\operatorname{WKB}}
\newcommand{\Pic}{\operatorname{Pic}}

\newcommand{\catc}{\mathsf{C}}

\newcommand{\catMod}{\mathsf{Mod}}
\newcommand{\catFun}[1][]{\mathsf{Hom}_{#1}}

\newcommand{\stkMod}[1][]{\mathfrak{Mod}_{#1}}
\newcommand{\stkAut}[1][]{\mathfrak{Aut}_{#1}}
\newcommand{\stkFun}[1][]{\mathfrak{Hom}_{#1}}

\newcommand{\stka}{\mathfrak{A}}

\newcommand{\stkg}{\mathfrak{G}}
\newcommand{\stkw}{\mathfrak{W}}
\newcommand{\stks}{\mathfrak{S}}

\newcommand{\astk}[1]{#1^+}

\newcommand{\approxto}[1][]{\xrightarrow[#1]{\approx}}

\newcommand{\W}{\shw^{\sqrt v}}

\begin{document}

\title[Classification of deformation quantization
algebroids\dots]{Classification of deformation quantization algebroids 
on complex symplectic manifolds}

\author[P. Polesello]{Pietro Polesello}
\address{Universit{\`a} di Padova\\
Dipartimento di Matematica Pura ed Applicata\\
via G. Belzoni, 7\\ 
35131 Padova, Italy}
\email{pietro@math.unipd.it}

\thanks{The author had the occasion of visiting Keio University and RIMS of 
Kyoto University during the preparation of this paper. Their hospitality is 
gratefully acknowledged .}

\date{}

\maketitle

\begin{abstract}
Deformation quantization algebroids over a complex symplectic manifold $X$ 
are locally given by rings of WKB operators, that is, microdifferential 
operators with an extra central parameter $\tau$. In this paper, we will 
show that such algebroids are classified by $H^2(X;k^*_X)$, where 
$k^*$ is a subgroup of the group of invertible formal power series in 
$\opb\tau$.\\

\noindent
{\bf Mathematics Subject Classification}: 46L65, 35A27, 18G5
\end{abstract}

\section*{Introduction} 
Let $M$ be a complex manifold, and $T^*M$ its cotangent bundle endowed with 
the canonical symplectic structure. Let $\shw_M$ be the sheaf of rings of WKB 
operators, that is, microdifferential operators with an extra central 
parameter $\tau$. This ring provides a  of $T^*M$. 
Recall that the order of the operators defines a filtration on $\shw_M$ such 
that its associated graded ring is isomorphic to $\O_{T^*M}[\tau^{-1},\tau]$.
Then, any filtered sheaf of rings which has $\O_{T^*M}[\tau^{-1},\tau]$ as 
graduate ring and which is locally isomorphic to $\shw_M$ gives another 
deformation quantization of $T^*M$. We call such an object a WKB-algebra.

On a complex symplectic manifold $X$ there may not exist a sheaf of rings of 
WKB operators, that is, a sheaf locally isomorphic to $\opb i \shw_M$, for any
symplectic local chart $i\colon X\supset U \to T^*M$. However, it is always 
defined an algebroid $\stkw_X$, which consists, roughly speaking, in 
considering the whole family of locally defined sheaves of WKB operators. 
This gives a  deformation quantization of $X$ 
(see \cite{Kashiwara1996,Kontsevich2001,Polesello-Schapira,D'Agnolo-Polesello2005}).
Again, the algebroid $\stkw_X$ is filtered and its associated graded is the trivial algebroid 
$\O_{X}[\tau^{-1},\tau]$. Then we may define a WKB-algebroid to be a filtered 
algebroid with the same graded as $\stkw_X$, and which is locally equivalent to 
$\stkw_X$. As before, any of these objects provides a deformation quantization of $X$.

The purpose of this paper is to show that WKB-algebroids are classified by 
$H^2(X;k^*_X)$, where $k^*$ is a subgroup of the group of invertible formal 
Laurent series in $\opb\tau$.

We refer to \cite{Deligne} for the classification of deformation quantization algebras 
on real symplectic manifolds.

The paper is organized as follows: we start by recalling the definition of 
WKB operator and that of WKB-algebra on $T^*M$, and by giving their 
classification. We then recall the main definitions and properties of filtered 
and graded stacks, and those of cohomology with values in a stack. With 
these tools at hand, we may define the WKB-algebroids on $X$ and give their 
classification.

\medskip
\noindent
{\bf Acknowledgement}
We wish to thank Masaki Kashiwara for useful suggestions.

\section{WKB-algebras}

The relation between Sato's microdifferential operators and WKB 
operators\footnote{WKB stands for Wentzel-Kramer-Brillouin.} is classical, and 
is discussed  e.g.~ in~\cite{Pham,AKKT}. We follow here the presentation 
in~\cite{Polesello-Schapira}, and we refer to~\cite{S-K-K,Kashiwara1986,Kashiwara2000} 
for the theory of microdifferential operators.

\medskip

Let $M$ be a complex manifold, and denote by $\rho\colon J^1M \to T^*M$ the 
projection from the 1-jet bundle to the cotangent bundle. 
Let $(t;\tau)$ be the system of homogeneous symplectic coordinates on $T^*\C$, and recall that $J^1 M$ is identified with the affine chart of the projective cotangent bundle $P^*(M\times \C)$ given by 
$\tau\neq 0$.
Denote by $\she_{M\times\C}$ the sheaf of microdifferential operators on $P^*(M\times\C)$.  
In a local coordinate system $(x,t)$ on $M\times\C$, consider the subring $\she_{M\times\C,\hat t}^{\sqrt v}$ of operators commuting with $\partial_t$. 
The ring of WKB operators is defined by
$$
\shw_M = \oim\rho (\she_{M\times\C,\hat t} |_{J^1M } ).
$$

In a local coordinate system $(x)$ on $M$, with associated symplectic local 
coordinates $(x;u)$ on $T^*M$, a WKB operator $P$ of order $m$
defined on a open subset $U$ of $T^*M$ has a total symbol
$$
\sigma(P)=\sum_{j=-\infty}^m p_j(x;u)\tau^{j},
$$
where the $p_j$'s are holomorphic functions on $U$ subject to the estimates
\begin{equation}\label{eq:estmicrod}
     \left\{ \begin{array}{l}
     \mbox{for any compact subset $K$ of $U$ there exists a constant}\\
     \mbox{$C_K>0$ such that for all $j<0$,}
     \sup\limits_{K}\vert p_{j}\vert \leq C_K^{-j}(-j)!.
     \end{array}\right.
\end{equation}
The product structure on $\shw_M$ is given by the Leibniz formula not 
involving $\tau$-derivatives. If $Q$ is another WKB operator defined on $U$ 
of total symbol $\sigma(Q)$, then 
$$\sigma(P\circ Q)=\sum_{\alpha\in\N^n} \frac{\tau^{-\vert\alpha\vert}}
{\alpha !} \partial^{\alpha}_u\sigma(P)\partial^{\alpha}_x\sigma(Q).
$$

\begin{remark}
The ring $\shw_M$ is a deformation quantization of $T^*M$ in the following 
sense. Setting $\hbar=\opb \tau$, the sheaf of formal WKB operators (obtained 
by dropping the estimates \eqref{eq:estmicrod}) of degree less than or equal 
to 0 is locally isomorphic to $\O_{T^*M}[\![\hbar]\!]$ as $\C_{T^*M}$-modules 
(via the total symbol), and it is equipped with an unitary associative product 
(the Leibniz rule) which induces a star-product on $\O_{T^*M}[\![\hbar]\!]$.
\end{remark}

Recall that the center of $\shw_M$ is the constant sheaf $k_{T^*M}$ with stalk
the subfield $k = \shw_{\operatorname{pt}} \subset \C[\![\tau^{-1},\tau]$ of 
WKB operators over a point, {\em i.e.} series $\sum_{j}
a_j{\tau}^j$ which satisfie the estimate:
\begin{equation*}
     \left\{ \begin{array}{l}
     \mbox{there exists a constant $C>0$ such that }\\
     \mbox{for all $j<0$, }\vert a_{j}\vert \leq C^{-j}(-j)!.
     \end{array}\right.
\end{equation*}

The sheaf  $\shw_M$ is filtered (over $\Z$), and one denotes by
$\shw_M(m)$ the sheaf of operators of order less than or equal to $m$. 
We denote by 
$$
\sigma_m(\cdot)\colon \shw_M(m)\to \shw_M(m)/\shw_M(m-1)
\simeq \sho_{T^*M}\cdot\tau^m
$$ 
the symbol map of order $m$. This function does not depend on the local 
coordinate system on $X$. If $\sigma_m(P)$ is not
identically zero, then one says that $P$ has order $m$ and
$\sigma_m(P)$ is called the principal symbol of $P$. 
In particular, an element $P$ in $\shw_M$ is invertible if and only if
its principal symbol is nowhere vanishing.
Moreover, the principal symbol map induces an isomorphism of graded rings:
$$
\sigma\colon\gr(\shw_M)\isoto\sho_{T^*M}[\opb\tau,\tau].
$$

\medskip

Let $\Omega_M$ be the canonical sheaf on $M$, that is, the sheaf of forms of 
top degree. Recall that each locally defined volume form $\theta\in\Omega_M$
gives rise to a local isomorphism $*_{\theta}\colon\shw_M^{\op} \isoto \shw_M$, 
which sends an operator $P$ to its formal adjoint $P^{*_{\theta}}$ with 
respect to $\theta$. Twisting $\shw_M$ by $\Omega_M$, one then gets a globally 
defined isomorphism of rings
$$
\shw_M^{\op} \isoto \opb {\pi}\Omega_M\tens\shw_M \tens \opb{\pi}
\Omega_M^{\tens -1}
\qquad P\mapsto \theta\tens P^{*_{\theta}}\tens \theta^{\tens -1}, 
$$
which does not depend on the choice of the volume form. 
(Here $\pi\colon T^*M\to M$ denotes the natural projection, 
$\Omega_M^{\tens -1}$  the $\O_M$-dual of $\Omega_M$
and the tensor product is over $\opb{\pi}\O_M$.) 
This leads to replace the ring $\shw_M$ by its twisted version 
by half-forms\footnote{Recall 
that the sections of $\W_M$ are locally defined by $\theta^{\tens 1/2}\tens 
P\tens \theta^{\tens -1/2}$ for a volume form $\theta$ and an operator $P$, 
with the equivalence relation $\theta_1^{\tens 1/2}\tens P_1\tens 
\theta_1^{\tens -1/2} = \theta_2^{\tens 1/2}\tens P_2\tens 
\theta_2^{\tens -1/2}$ if and only if $P_2 = (\theta_1/\theta_2)^{1/2} P_1 
(\theta_1/\theta_2)^{-1/2}$.}
$$
\W_M =
\opb {\pi} \Omega_M^{\tens 1/2}\tens\shw_M\tens\opb{\pi}\Omega_M^{\tens -1/2}.
$$

The $k$-algebra $\W_M$ is locally isomorphic to $\shw_M$ and has the 
following properties:
\begin{itemize}
\item[(i)] it is filtered;
\item[(ii)] there is an isomorphism of graded rings
\begin{equation*}
\sigma\colon\gr(\W_M)\isoto\sho_{T^*M}[\opb\tau,\tau];
\end{equation*}
\item[(iii)] it is endowed with an anti-involution, {\em i.e.} an isomorphism 
of rings $$*\colon (\W_M)^{\op}\isoto\W_M \quad \mbox{such that 
$*^2=\id$.}$$
\end{itemize}

This suggests the following
\begin{definition}
A WKB-algebra on $T^*M$ is a sheaf of $k$-algebras $\sha$ together with
\begin{itemize}
\item[(i)] a filtration $\{F_m\sha\}_{m\in\Z}$; 
\item[(ii)] an isomorphism of graded rings $\nu\colon\gr(\sha)\isoto
\O_{T^*M}[\opb\tau,\tau]$; 
\item[(iii)] an anti-involution $\iota$;
\end{itemize} 
such that the triplet $(\sha,\nu,\iota)$ is locally isomorphic to 
$(\W_M, \sigma, *)$. 

A morphism of WKB-algebras is a $k$-algebra morphism compatible with the 
structures (i), (ii) and (iii).
\end{definition}

By definition, an isomorphism of WKB-algebras $\varphi\colon \sha_1\to \sha_2$ 
is a $k$-algebra isomorphism commuting with the anti-involutions, mapping 
$F_m\sha_1$ to $F_m\sha_2$ in such a way that $\nu^2_m(\varphi(P)) =
\nu^1_m(P)$ for all $P\in F_m\sha_1$. (Here  $\nu^i_m$ denotes the symbol map 
$F_m(\sha_i)\to F_m(\sha_i)/F_{m-1}(\sha_i)\simeq  
\sho_{T^*M}\cdot\tau^m$ of order $m$, for $i=1,2$.) 
This translates to WKB operators the notion of equivalence between star-products.
Hence any (formal) WKB-algebra provides a deformation quantization 
of $T^*M$. See \cite{BoutetdeMonvel1999,BoutetdeMonvel2002} for similar 
definitions in the context of microdifferential and Toeplitz operators. 

\begin{example}
Let $f\colon T^*M \to T^*M$ be a symplectic transformation. Then $\W_M$ induces an
anti-involution on $\opb f \W_M$ and a filtration such that the associated graded ring is 
isomorphic (via $f$) to $\O_{T^*M}[\opb\tau,\tau]$. By \cite{Polesello-Schapira}, locally there 
exists a Quantized Symplectic Transformation over $f$, that is, an isomorphism 
$\opb f \W_M\simeq \W_M$of filtered $k$-algebras, commuting with the anti-involutions
and which preserves the graded rings. It follows that $\opb f \W_M$ is a WKB-algebra.
\end{example}

Denote by $\shaut[\WKB](\W_M)$ the group of WKB-algebra automorphisms of 
$\W_M$ and set
\begin{equation*}
\begin{split}\shw^{\sqrt v, *}_M & = \{P\in\W_M; 
\mbox{ $P$ has order 0, $\sigma_0(P)=1$ and  $PP^*=1$}\},\\
k^* & = \{s(\tau)\in k;\mbox{ $s(\tau)=1+\sum_{j< 0}a_j{\tau}^j$ and 
$s(\tau)s(-\tau) =1$}\}.
\end{split}
\end{equation*}
Note that $\shw^{\sqrt v, *}_M$ is a subgroup of the group 
$\shw^{\sqrt v, \times}_M$ of invertible WKB operators, and that 
$k^*  = \shw^{\sqrt v, *}_{\operatorname{pt}}$.

\begin{lemma}(cf \cite{Polesello-Schapira})\label{lemma:key}
There is an exact sequence of groups on 
$T^*M$
\begin{equation}\label{eq:key}
1\to k^*_{T^*M} \to \shw^{\sqrt v, *}_M \to[\ad]  \shaut[\WKB](\W_M) \to 1,
\end{equation}
where $\ad(P)(Q)=PQP^{-1}$ for any $P\in \shw^{\sqrt v, *}_M$ and  
$Q\in\W_M$.
\end{lemma}


The set of isomorphism classes of WKB-algebras on $T^*M$ is in bijection with 
$H^1(T^*M; \shaut[\WKB](\W_M))$. Hence we get
\begin{corollary}
WKB-algebras on $T^*M$ are classified by the pointed set
$H^1(T^*M;\shw^{\sqrt v, *}_M/k^*_{T^*M})$.
\end{corollary}

\section{Filtered and graded stacks}

To define WKB-algebroids, we need to translate the notions of filtration and 
graduation from sheaves to stacks. We start here by recalling what a filtered 
(resp. graded) category is and how to associate a graded category to a 
filtered one. Then we stakify these definitions.
We assume that the reader is familiar with the basic notions from the theory 
of stacks which are, roughly speaking, sheaves of categories. (The classical 
reference is~\cite{Giraud1971}, and a short presentation is given {\em e.g.}
in~\cite{Kashiwara1996,D'Agnolo-Polesello2003}.)

\medskip
Let $R$ be a commutative ring. 
\begin{definition}\label{def:filtered} 
     A filtered (resp. graded) $R$-category is an 
     $R$-category\footnote{Recall that an $R$-category is a category whose 
     sets of morphisms are endowed with an $R$-module structure, so that 
     composition is bilinear. An $R$-functor is a functor between 
     $R$-categories which is linear at the level of morphisms.} 
     $\catc$ such that:
     \begin{itemize}
     \item[$\bullet$] for any objects $P,Q\in\catc$,  the $R$-module 
     $\Hom[\catc](P,Q)$ is filtered (resp. graded) over $\Z$;
     \item[$\bullet$] for any $P,Q,R\in \catc$ and any morphisms $f$ in 
     $F_m\Hom[\catc](Q,R)$ (resp. in $G_m\Hom[\catc](Q,R)$) and $g$ in 
     $F_n\Hom[\catc](P,Q)$ (resp. in $G_n\Hom[\catc](P,Q)$), the composition 
     $f\circ g$ is in $F_{m+n}\Hom[\catc](P,R)$ (resp. in 
     $G_{m+n}\Hom[\catc](P,R)$);
     \item[$\bullet$] for each $P\in \catc$, the identity morphism $\id_{P}$ 
     is in $F_0\Hom[\catc](P,P)$ (resp. in $G_0\Hom[\catc](P,P)$). 
     \end{itemize}
     A filtered (resp. graded) $R$-functor is an $R$-functor which respects 
     the filtrations (resp. graduations) at the level of morphims.
\end{definition}

To any filtered $R$-category $\catc$ there is an associated graded $R$-category 
$\Gr(\catc)$, whose objects are the same of those of $\catc$ and the
morphisms are defined by $\Hom[\Gr(\catc)](P,Q) = \Gr(\Hom[\catc](P,Q))$ for any objects 
$P,Q$. In this way, we get a functor from filtered $R$-categories to graded ones.

Following the presentation in~\cite{D'Agnolo-Polesello2005}, recall that there 
is a fully faithful functor from filtered (resp. graded) 
$R$-algebras to filtered (resp. graded) $R$-categories, which sends a filtered 
(resp. graded) $R$-algebra $A$ the category $\astk A$ with a single 
object $\bullet$  and $\Endo(\bullet)=A$ as set of morphisms. Hence, the 
functors $\Gr$ and $\astk{}$ commutes, that is, for any filtered $R$-algebra 
$A$ one has $\Gr (\astk A)=\astk{\Gr (A)}$.

Note that, if $A$ is a filtered $R$-algebra, the $R$-category $\catMod_F(A)$ 
of filtered left $A$-modules has a natural filtration: for any filtered 
$A$-modules $M$ and $N$, one sets $F_m\Hom[FA](M,N)=\Hom[FA](M,N(m))$,
where $N(m)$ has the same underlying $A$-module as $N$, and the filtration 
is given by $F_n N(m) = F_{n+m}N$. 
One easily checks that $\catMod_F(A)$ is equivalent to the category 
$\catFun[F](\astk A,\catMod_F(R))$ of filtered $R$-functors from 
$\astk A$ to $\catMod_F(R)$ and that the Yoneda embedding
$$
\astk A \to \catFun[F]((\astk A)^\op, \catMod_F(R)) \approx 
\catMod_F(A^\op)
$$
identifies $\astk A$ with the full subcategory of filtered right $A$-modules 
which are free of rank one. Everything remains true replacing filtered 
algebras and categories by graded ones.


\medskip

Let $X$ be a topological space, and $\shr$ a sheaf of commutative rings.

As for categories, there are natural notions of filtered (resp. graded) 
$\shr$-stack, and of filtered (resp. graded) $\shr$-functor between 
filtered (resp.  graded) $\shr$-stacks.

 
As above, we denote by $\astk{}$ the (faithful and locally full) functor from filtered (resp. graded) 
$\shr$-algebras to filtered (resp. graded) $\shr$-categories, which sends
a filtered (resp. graded) $\shr$-algebra $\sha$ to the stack $\astk\sha$ 
defined as follows: it is the stack associated with the separated prestack 
$X \supset U\mapsto \astk{\sha(U)}$.


If $\sha$ is a filtered $\shr$-algebra, then the stack 
$\stkMod[F](\sha)$ of filtered left $\sha$-modules is filtered and equivalent 
to the stack of filtered functors $\stkFun[F](\astk \sha,\stkMod[F](\shr))$, 
and the Yoneda embedding gives a fully faithful functor 
\begin{equation*}
\label{eq:Yoneda}
\astk \sha \to \stkFun[F]((\astk \sha)^\op,\stkMod[F](\shr)) \approx
\stkMod[F](\sha^\op)
\end{equation*}
into the stack of filtered right $\sha$-modules. This identifies $\astk\sha$ 
with the full substack of locally free filtered right $\sha$-modules of rank 
one. As above, everything remains true replacing filtered algebras and stacks 
by graded ones.

Let $\stks$ be a filtered $\shr$-stack. We denote by $\gr (\stks)$ the graded stack 
associated to the pre-stack $X \supset U\mapsto \Gr (\stks(U))$.
\begin{proposition}
Let $\sha$ be a filtered $\shr$-algebra and $\gr (\sha)$ its associated graded 
ring. Then there is an equivalence of graded 
stacks $\gr (\sha^+)\approx\gr (\sha)^+$. 
\end{proposition}

\begin{proof}
Let $\shl$ be a locally free right filtered $\sha$-module of rank one (that is,
an object of $\sha^+$). Its associated graded module $\gr (\shl)$ is a 
locally free right graded $\gr (\sha)$-module of rank one (that is, an object 
of $\gr( \sha)^+$). Hence the assignement $\shl \mapsto \\gr (\shl)$ induces  
a functor $\gr( \sha^+)\to\gr( \sha)^+$ of graded stacks.
Since at each $x\in X$ this reduces to the equality $\Gr (\astk {\sha_x}) = 
\astk{\Gr (\sha_x)}$, it follows that it is a global equivalence. 
\end{proof}

Recall from \cite{Kontsevich2001,D'Agnolo-Polesello2005} that an 
$\shr$-algebroid stack is an $\shr$-stack $\stka$ which is locally non-empty and 
locally connected by isomorphisms. Equivalently, for any  
$x\in X$ there exist an open subset $U\subset X$ containing $x$ and an 
$\shr$-algebra $\sha$ on $U$ such that $\stka|_U\approx \astk \sha$. 
\begin{corollary}
Let $\stka$ be a filtered $\shr$-stack.  If $\stka$ is an $\shr$-algebroid stack, then it associated 
graded stack $\gr(\stka)$ is again an $\shr$-algebroid stack.
\end{corollary}

\section{WKB-algebroids}

Let $(X,\omega)$ be a complex symplectic manifold. Recall that a local model 
for $X$ is an open subset $U$ of the cotangent bundle $T^*M$ of a complex 
manifold $M$, equipped with the canonical symplectic form. 
Although there may not exist a globally defined WKB-algebra on $X$, that is, 
a sheaf locally isomorphic to $\opb i \shw_M$ for any symplectic local chart 
$i\colon X\supset U \to T^*M$, Polesello-Schapira~\cite{Polesello-Schapira}
defined a canonical stack of WKB-modules on $X$. 
Following~\cite{D'Agnolo-Polesello2005}, this result may be restated as:

\begin{theorem}
On any complex symplectic manifold $X$ there exists a canonical $k$-stack 
$\stkw_X$ which is locally equivalent to $\astk{(\opb i \W_M)}$ 
for any symplectic local chart $i\colon X\supset U \to T^*M$.
\end{theorem}

By definition, $\stkw_X$ is a $k$-algebroid stack. Hence there exists a canonical
WKB-algebra on $X$ if and only if $\stkw_X$ has a global object. 

\begin{proposition}
The $k$-algebroid stack $\stkw_X$ has the following properties:
\begin{itemize}
\item[(i')] it is filtered;
\item[(ii')] there is a natural equivalence of graded stacks 
$$\pmb\sigma\colon\gr(\stkw_X)\approxto \astk {(\O_X[\tau^{-1},\tau])};$$
\item[(iii')] it is endowed with an anti-involution $*$, that is, with a linear 
equivalence 
$$\pmb*\colon \stkw_X^\op\approxto \stkw_X$$
and an invertible transformation $\epsilon\colon \pmb*^2 \Rightarrow\id_{\stkw_X}$ 
such that the transformations $\epsilon\id_{\pmb*}\colon \pmb*^3 \Rightarrow \pmb*$ and 
$\id_{\pmb*}\epsilon\colon \pmb* \Rightarrow \pmb*^3$ are inverse one to each other.
\end{itemize} 
\end{proposition}

We may mimic the definition of WKB-algebra and get the following
\begin{definition}
A WKB-algebroid on $X$ is a $k$-stack $\stka$ endowed with 
\begin{itemize}
\item[$\astk {(i)}$] a filtration;
\item[$\astk {(ii)}$] an equivalence of graded stacks 
$\pmb\nu\colon\gr(\stka)\approxto \astk {(\O_X[\tau^{-1},\tau])};$
\item[$\astk {(iii)}$] an anti-involution $\pmb\iota$;
\end{itemize} 
such that the triplet $(\stka,\pmb\nu,\pmb\iota)$ is locally equivalent to 
$(\stkw_X,\pmb\sigma,\pmb*)$.

A functor of WKB-algebroids is a $k$-functor compatible with the 
structures $\astk {(i)}$, $\astk {(ii)}$ and $\astk {(iii)}$.
\end{definition}
As (formal) WKB-algebras give the deformation quantizations of $T^*M$, we may 
say that (formal) WKB-algebroids provide the deformation quantizations of $X$.

\begin{definition}
We call $\stkw_X$ the canonical WKB-algebroid on $X$.
\end{definition}

\section{Cohomology with values in a stack}

As for classifying WKB-algebras one uses cohomology with values in a sheaf of 
groups, so to classify WKB-algebroids we need a cohomology theory with values 
in a stack with group-like properties. In this section we briefly recall the 
definition of cohomology with values in a stack and show how to describe it 
explicitly by means of the notion of crossed module.
References are made to \cite{Breen1992,Breen1994}.

We assume that the reader is familiar with the notions of monoidal category 
and monoidal functor. (The classical reference is \cite{MacLane}.)

Let $X$ be a topological space.

\begin{definition}
\begin{itemize}
     \item[(i)] A 2-group\footnote{We follow here the terminology of 
     Baez-Lauda [{\em Higher-dimensional algebra V: 2-groups}, e-print (2004) 
     \texttt{arXiv:math.QA/0307200}], which seems to us more 
     friendly than the classical one of $gr$-category due to Grothendieck.} 
     is a rigid monoidal groupoid, {\em i.e.} a monoidal category $(\mathsf{G},
     \tens,{\bf 1})$ with all the morphisms invertible and such that for any 
     object $P\in \mathsf{G}$ there exist an object $Q$ and natural 
     morphisms $P\tens Q\simeq {\bf 1}$ and $Q\tens P\simeq {\bf 1}$. 
     A functor of 2-groups is a monoidal functor between the underlying 
     monoidal categories.
     
     \item[(ii)] A pre-stack (resp. stack) of 2-groups on $X$ is a pre-stack (resp. stack) 
     $\stkg$ such that for each open subset $U\subset X$, the category $\stkg(U)$ is 
     2-group and the restriction functors are functors of 2-groups.
    \end{itemize}
\end{definition}
If there is no risk of confusion, a stack of 2-groups on $X$ will be simply called a 2-group on $X$.

\begin{example}
Let $\shg$ be a sheaf of groups on $X$.
\begin{itemize}
     \item[(i)] The discrete stack $\shg[0]$ defined by trivially enriching  
     $\shg$ with identity arrows is a 2-group on $X$.
     
     \item[(ii)]  Let $\shg[1]$ be the stack in groupoids associated to the 
     separated pre-stack whose category on an open subset $U\subset X$ 
     has a single object $\bullet$  and $\Endo(\bullet)=\shg(U)$ as set of 
     morphisms. Then $\shg[1]$ is equivalent to the stack of right 
     $\shg$-torsors and it defines a 2-group on $X$ if and only if $\shg$ is 
     commutative.
     
\end{itemize}
\end{example}

Let $\stkg$ be a pre-stack of 2-groups on $X$. 
We define the 0-th cohomology group of $X$ with values in $\stkg$ to be 
$$ H^0(X;\stkg)=\ilim[\shu] H^0(\shu;\stkg),$$ 
where $\shu$ ranges over open coverings of $X$. 
For an open covering $\shu = \{U_{i}\}_{i\in I}$, the elements of 
$H^0(\shu;\stkg)$ are represented by pairs $(\{\shp_{i}\}, \{\alpha_{ij}\})$ 
(the 0-cocycles), where $\shp_{i}$ is an object in $\stkg (U_{i})$ and 
$\alpha_{ij}\colon \shp_{j} \isoto \shp_{i}$ is an isomorphism on double 
intersection $U_{ij} = U_i\cap U_j$, such that 
$\alpha_{ij}\circ \alpha_{jk}=\alpha_{ik}$ on triple intersection $U_{ijk}$, 
with the relation $(\{\shp_{i}\}, \{\alpha_{ij}\})$ is equivalent to 
$(\{\shp'_{i}\}, \{\alpha'_{ij}\})$ 
if and only if there exists an isomorphism 
$\delta_{i}\colon \shp'_{i} \isoto \shp_{i}$ compatible with 
$\alpha_{ij}$ and $\alpha'_{ij}$ on $U_{ij}$. 

Note that, if $\stkg$ is a stack of 2-groups, then $H^0(X;\stkg)$ is 
isomorphic to the group of isomorphism classes of objects in $\stkg(X)$.

\medskip

Similarly, the 1-st  cohomology (pointed) set of $X$ with values in 
$\stkg$ is defined as
$$H^1(X;\stkg)=\ilim[\shu] H^1(\shu;\stkg),$$
where $\shu$ ranges over open coverings of $X$.
For an open covering $\shu = \{U_{i}\}_{i\in I}$, the elements of 
$H^1(\shu;\stkg)$ are given by pairs $(\{\shp_{ij}\}, \{\alpha_{ijk}\})$ 
(the 1-cocycles), where $\shp_{ij}$ is an object in $\stkg (U_{ij})$ and 
$\alpha_{ijk}\colon \shp_{ij} \tens \shp_{jk} \isoto \shp_{ik}$ is an 
isomorphism on $U_{ijk}$ such that the diagram on quadruple intersection 
$U_{ijkl}$
\begin{equation*}
\xymatrix@C5em{ \shp_{ij}\tens \shp_{jk}\tens \shp_{kl} 
\ar[r]^-{\alpha_{ijk}\tens\id_{\shp_{kl}}} \ar[d]^{\id_{\shp_{ij}}\tens\alpha_{jkl}} & 
\shp_{ik}\tens \shp_{kl}\ar[d]^{\alpha_{ikl}} \\
\shp_{ij}\tens \shp_{jl} \ar[r]^-{\alpha_{ijl}} 
& \shp_{il} }
\end{equation*}
commutes. The 1-cocylces $(\{\shp_{ij}\}, \{\alpha_{ijk}\})$ and 
$(\{\shp'_{ij}\}, \{\alpha'_{ijk}\})$ are equivalent if and only if there 
exists a pair $(\{\shq_i\}, \{\delta_{ij}\})$, with $\shq_i$ an object of 
$\stkg(U_i)$ and $\delta_{ij}\colon \shp'_{ij} \tens\shq_j \isoto 
\shq_i \tens \shp_{ij}$ an isomorphism on $U_{ij}$ such that the diagram on 
$U_{ijk}$
\begin{equation*}
\xymatrix@C4em@R3em{ \shp'_{ij}\tens \shp'_{jk}\tens \shq_k
\ar[r]^-{\id_{\shp'_{ij}}\tens\delta_{jk}} 
\ar[d]^{\alpha'_{ijk}\tens\id_{\shq_{k}}} & 
\shp'_{ij}\tens\shq_j\tens\shp_{jk} 
\ar[r]^{\delta_{ij}\tens\id_{\shp_{jk}}} & 
\shq_i\tens \shp_{ij}\tens \shp_{jk} 
\ar[d]^{\id_{\shq_{i}}\tens\alpha_{ijk}} \\
\shp'_{ik}\tens \shq_k  \ar[rr]^-{\delta_{ik}} 
&& \shq_i\tens \shp_{ik} }
\end{equation*}
commutes. 

\medskip

In the rest of the section we will give a more explicit description of the cohomology 
with values in a stack by means of cocycles with values in a crossed module. 
(This was Breen's approach to non abelian cohomology of Giraud~\cite{Giraud1971}.)

\begin{definition}
A crossed module on $X$ is a complex of sheaves of groups 
$\shg^{-1}\to[d]\shg^0$ endowed with a left action of $\shg^0$ on $\shg^{-1}$
such that for any local sections $g\in\shg^0$ and $h,h'\in\shg^{-1}$ one has
$$
d({}^gh)=\ad(g)(d(h)) \qquad {}^{d(h')} h=\ad(h')(h).
$$
(Here we use the convention as in \cite{Breen1994}
for which $\shg^{i}$ is in $i$-th degree.)

A morphism of crossed modules is a morphism of complexes compatible with the 
actions in the natural way. 
\end{definition}

Associated to each crossed module $\shg^{-1}\to[d]\shg^0$ there is 
2-group on $X$, which we denote by $[\shg^{-1}\to[d]\shg^0]$, 
defined as follows: it is the stack associated to the separated 
pre-stack of 2-groups whose objects on an open subset $U\subset X$ are the
sections $g\in \shg^0(U)$ with 2-group law $g\tens g'=gg'$, and whose 
morphisms $g\to g'$ are given by sections $h\in \shg^{-1}(U)$ such that 
$g' = d(h) g$,  with the 2-group structure given by $(g_1\to[h_1] g'_1)\tens 
(g_2\to[h_2] g'_2) = g_1g_2\to[h_1{}^{g_1}h_2] g'_1g'_2$.
Similarly, each morphism of crossed modules induces a functor of the 
corresponding 2-groups.

\begin{remark}
In fact, it is true that any 2-group on $X$ comes from a crossed module.
However, this result is not of practical use. We refer to 
\cite{SGA4} for the proof of this fact in the commutative case and to 
\cite{Brown-Spencer} for the non commutative case on $X=\operatorname{pt}$. 
\end{remark}

%

\begin{example}
Let $\shg$ be a sheaf of groups on $X$.
\begin{itemize}
     \item[(i)]The 2-group defined by the 
     crossed module $1 \to \shg$ is identified with $\shg[0]$. 
     
     \item[(ii)]  If moreover $\shg$ is commutative, the complex $\shg \to 1$ 
     is a crossed module and its associated 2-group is identified with 
     $\shg[1]$. 
     
\end{itemize}    
\end{example}

Let $\shg^{-1}\to[d]\shg^0$ be a crossed module on $X$. Then the cohomology 
of $X$ with values in the 2-group $[\shg^{-1}\to[d]\shg^0]$ admits a very 
explicit description, which we recall below. This is usually referred as the 
(hyper-)cohomology of $X$ with values in $\shg^{-1}\to[d]\shg^0$.

By definition, an object $\shp$ of $[\shg^{-1}\to[d]\shg^0]$ on 
an open subset $U\subset X$ is described by an open covering $U = 
\Union\nolimits_i U_i$ and sections $\{g_i\}\in \shg^0(U_i)$, subject 
to the relation $g_i = d(h_{ij}) g_j$ on double intersections $U_{ij}$, for 
given sections $\{h_{ij}\}\in \shg^{-1}(U_{ij})$ satisfying $h_{ij}h_{jk}=h_{ik}$ 
on triple intersections $U_{ijk}$.
Hence, up to a refinement of the open covering $\shu=\{U_{i}\}_{i\in I}$ of 
$X$, the 0-cocycles on $\shu$ with values in $[\shg^{-1}\to[d]\shg^0]$, may be 
described by pairs $(\{g_i\},\{h_{ij}\})$, where $g_i\in \shg^0(U_i)$ 
and $h_{ij}\in \shg^{-1}(U_{ij})$ are sections satisfying the relations
\begin{equation*}
\label{nonab1}
\begin{cases}
g_i=d(h_{ij})g_j \quad \text{in }\shg^0(U_{ij})\\
h_{ij}h_{jk} = h_{ik} \quad \text{ in } \shg^{-1}(U_{ijk}),
\end{cases}
\end{equation*}
and $(\{g_i\},\{h_{ij}\})$ is equivalent to $(\{g'_i\},\{h'_{ij}\})$ if 
and only if there exist sections $\{k_i\}\in\shg^{-1}(U_i)$ such that the
following relations hold
\begin{equation*}
\label{nonabeq1}
\begin{cases}
g'_i = d(k_i)g_i  \\
h'_{ij}  k_j = k_i h_{ij} .
\end{cases}
\end{equation*}

\medskip

The same description for 1-cocycles needs some care, since one has to 
consider open coverings for any double intersection $U_{ij}$. In other words, 
one has to replace coverings by hypercoverings. Indices become thus 
very cumbersome, and we will not write them explicitly\footnote{Recall that, 
on a paracompact space, usual coverings are cofinal among hypercoverings}.
Hence the 1-cocycles on $\shu$ with values in $[\shg^{-1}\to[d]\shg^0]$, may 
be described by pairs $(\{g_{ij}\},\{h_{ijk}\})$, with $g_{ij}\in 
\shg^0(U_{ij})$ and $h_{ijk}\in \shg^{-1}(U_{ijk})$ satisfying the
relations
\begin{equation*}
\label{nonab2}
\begin{cases}
g_{ij}g_{jk}=d(h_{ijk})g_{ik} \quad \text{  in } \shg^0(U_{ijk})\\
h_{ijk}h_{ikl} = {}^{g_{ij}}h_{jkl}h_{ijl} \quad \text{in } 
\shg^{-1}(U_{ijkl}).
\end{cases}
\end{equation*}
Moreover, $(\{g_{ij}\},\{h_{ijk}\})$ is equivalent to 
$(\{g'_{ij}\},\{h'_{ijk}\})$ if and only if there exists a pair 
$(\{l_i\},\{k_{ij}\})$, with $k_{ij}\in\shg^{-1}(U_{ij})$ and 
$l_i\in\shg^0(U_i)$, such that 
\begin{equation*}
\label{nonabeq2}
\begin{cases}
g'_{ij}  l_j = d(k_{ij})l_i g_{ij} \\
h'_{ijk} k_{ik} = {}^{g'_{ij}} k_{jk}k_{ij} {}^{l_i} h_{ijk}.
\end{cases}
\end{equation*}

\medskip

Taking the extremal cases $[1\to\shg]$ and $[\shg\to 1]$, the latter when the 
group $\shg$ is commutative, one easily recovers from the previous description 
the usual definition of the Cech cohomology of $X$ with values in $\shg$. 
Hence one has the following
\begin{proposition}\label{prop:hyper}
Let $\shg$ be a sheaf of groups on $X$. Then there is an isomorphism 
(of groups if $i=0$, of pointed sets if $i=1$) 
$$
H^i(X;\shg[0])\simeq  H^{i}(X;\shg).
$$
If moreover $\shg$ is commutative, then there are isomorphisms of groups (for 
$i=0,1$)
$$
H^i(X;\shg[1])\simeq  H^{i+1}(X;\shg).
$$
\end{proposition}

\section{Classification of WKB-algebroids}

Let $(X,\omega)$ be a complex symplectic manifold of dimension $2n$ and $\stkw_X$ the 
canonical WKB-algebroid on $X$.

Let $\stka$ be another WKB-algebroid. By definition, there exists an open 
covering $X=\Union\nolimits_i U_i$ such that $\stka|_{U_i}$ is equivalent to 
$\stkw_X|_{U_i}$ as WKB-algebroids.
Let $\Phi_i \colon \stka|_{U_i} \to \stkw_X|_{U_i}$ and 
$\Psi_i \colon \stkw_X|_{U_i}\to \stka|_{U_i}$ be quasi-inverse to each other.
On double intersections $U_{ij}$ there are WKB-algebroid 
equivalences $\Phi_{ij} = \Phi_i\Psi_j \colon \stkw_X|_{U_{ij}} \to 
\stkw_X|_{U_{ij}}$, and on triple intersections $U_{ijk}$ there are 
invertible transformations $\alpha_{ijk} \colon \Phi_{ij}\Phi_{jk}
\Rightarrow \Phi_{ik}$ induced by $\Psi_j\Phi_j\Rightarrow\id$. 
On quadruple intersections $U_{ijkl}$ the following diagram commutes
\begin{equation}
\label{eq:alpha}
\xymatrix@C5em{
\Phi_{ij}\Phi_{jk}\Phi_{kl} \ar@{=>}[r]^{\alpha_{ijk}\id_{\Phi_{kl}}} \ar@{=>}[d]^{\id_{\Phi_{ij}}\alpha_{jkl}}
& \Phi_{ik}\Phi_{kl} \ar@{=>}[d]^{\alpha_{ikl}} \\
\Phi_{ij}\Phi_{jl} \ar@{=>}[r]^{\alpha_{ijl}} & \Phi_{il} .
}
\end{equation}

It follows that WKB-algebroids are described by 1-cocycles 
$(\Phi_{ij},\alpha_{ijk})$ with values in the stack of 2-groups 
$\stkAut[\WKB](\stkw_X)^\times$ of autoequivalences of $\stkw_X$ as 
WKB-algebroid. (Here the upper index $\times$ means that all the non-invertible 
morphisms have been removed.)

Denote by $\WKB(X)$ the set of equivalence classes of WKB-algebroid 
on $X$, pointed by the class of $\stkw_X$. Hence one gets an isomorphism of
pointed sets 
$$
\WKB(X)\simeq H^1(X;\stkAut[\WKB](\stkw_X)^{\times}).
$$ 

\medskip

Let us briefly recall how to describe more explicitly the 1-cocycle $(\Phi_{ij},\alpha_{ijk})$ 
attached to a WKB-algebroid $\stka$. We follow \cite{Polesello-Schapira,D'Agnolo-Polesello2005}.

By definition, $\stkw_X$ is locally equivalent to 
$\astk{(\opb f \W_M)}$ for any symplectic local chart 
$f\colon X\supset U \to T^*M$. Hence, up to a refinement of the open covering 
$X=\Union\nolimits_i U_i$, one may suppose that $\stkw_X$ is equivalent on 
$U_i$ to $\astk {\W_i{}} = \astk{(\opb{f_i}\W_M)}$, for 
a symplectic embedding $f_i \colon U_i \to T^*M$ with $M=\C^n$. 
On $U_{ij}$ the functor $\Phi_{ij}\colon \astk{\W_j{}} \to 
\astk{\W_i{}}$ is then locally induced by WKB-algebra isomorphisms. 
Shrinking again the open covering, we may find an isomorphism of WKB-algebras 
$\varphi_{ij} \colon \W_j \to \W_i$ on $U_{ij}$ such 
that $\astk{\varphi_{ij}} = \Phi_{ij}|_{U_{ij}}$. 
On $U_{ijk}$ we have an invertible transformation 
$\alpha_{ijk} \colon \astk{\varphi_{ij}}\astk{\varphi_{jk}} \Rightarrow 
\astk{\varphi_{ik}}$, so that there exist a section 
$P_{ijk} \in \shw^{\sqrt v, *}_i$ such that
$$
\varphi_{ij}\varphi_{jk}= \ad(P_{ijk}) \varphi_{ik}.
$$
Finally, on $U_{ijkl}$ the diagram \eqref{eq:alpha} corresponds to the 
equality
$$
P_{ijk} P_{ikl} = \varphi_{ij}(P_{jkl}) P_{ijl}.
$$

The datum of $(\{f_i\},\{\varphi_{ij}\}, \{P_{ijk}\})$ as above is enough 
to reconstruct $\stka$ (up to equivalence). 


\medskip

In the particular case of $X=T^*M$, one has $f_i=\id$. A direct computation 
as above shows that there is  an equivalence of 2-groups
\begin{equation}\label{WKB-aut}
\stkAut[\WKB](\astk{(\W_M)})^\times
\approx \left[\shw^{\sqrt v, *}_M \to[\ad]  \shaut[\WKB](\W_M) \right].
\end{equation}

\bigskip

We are now ready to prove the following

\begin{theorem}\label{th:classification}
There is an isomorphism of pointed sets
$$\WKB(X)\simeq H^2(X;k^*_X).$$
\end{theorem}

\begin{proof}
Consider the natural functor of 2-groups
$$
F\colon k^*_X[1] \longrightarrow\stkAut[\WKB](\stkw_X)^{\times}
$$
induced by the functor of pre-stacks which sends the unique object 
$\bullet $ to the identity functor $\id_{\stkw_X}$.
At any point $p\in X$, we may find a symplectic local chart 
$i\colon X\supset U \to T^*M$ around $p$, such that $\stkw_X|_U$ is equivalent 
to $ \astk{\W_U{}}$ as WKB-algebroid. (Here we set $\W_U =\opb i \W_M$.)
We thus have a chain of equivalences of 2-groups
\begin{equation*}
\begin{split}
\stkAut[\WKB](\stkw_X)^\times|_U 
& \approx  \stkAut[\WKB](\stkw_X|_U)^\times\\
& \approx \stkAut[\WKB](\astk{(\shw_U^{\sqrt v})})^\times\\
& \approx \left[ \shw^{\sqrt v, *}_U \to[\ad]  \shaut[\WKB](\W_U) \right], \\
\end{split}
\end{equation*}
(the last one follows from~\eqref{WKB-aut}) and hence the functor $F$ 
restricts on $U$ to
$$
F|_U\colon k^*_U[1] \longrightarrow \left[\shw^{\sqrt v, *}_U\to[\ad] 
\shaut[\WKB](\W_U) \right].
$$
By Lemma \ref{lemma:key2}, this is an equivalence of 2-groups so that the 
functor $F$ is locally, and hence globally, an equivalence.
We thus get a chain of isomorphisms of pointed sets
\begin{equation*}
H^1(X;\stkAut[\WKB](\stkw_X)^{\times}) \simeq H^1(X;k^*_X[1]) 
\simeq  H^2(X;k^*_X),
\end{equation*}
where the latter follows by Proposition \ref{prop:hyper}. 
\end{proof}

\begin{lemma}\label{lemma:key2}
Let $M$ be a complex manifold. Then there is an equivalence of 2-groups on 
$T^*M$
$$
k^*_{T^*M}[1] \approx \left[\shw_M^{\sqrt v, *} \to[\ad] 
\shaut[\WKB](\W_M) \right].
$$
\end{lemma}

\begin{proof}
This follows by a direct computation, using \eqref{eq:key}.
\end{proof}

\begin{corollary}
The set $\WKB({X})$ has an abelian group structure.
\end{corollary}

Note that, by much the same proof of Theorem~\ref{th:classification}, one gets 
an isomorphism between the group $\Pic_{\WKB}(\stkw_{X})$ of isomorphism 
classes of autoequivalences\footnote{For a sheaf of rings $\sha$, one usually
denotes by $\Pic(\sha)$ the group of isomorphism classes of invertible 
$\sha\tens[\shr]\sha^{\op}$-modules. This is consistent with our notation, 
since by Morita theorem $\Pic(\sha)$ is isomorphic to the group of isomorphism 
classes of linear autoequivalences of $\stkMod(\sha)$.} of $\stkw_{X}$ as 
WKB-algebroid and $H^1(X;k^*_{X})$.

\begin{remark}
If $X=T^*M$ for a complex manifold $M$, the coboundary map 
$$
\delta\colon H^1(T^*M;\shw^{\sqrt v, *}_M/k^*_{T^*M}) \to H^2(T^*M;k^*_{T^*M})
$$
associated to the exact sequence~\eqref{eq:key}, may be interpreted as the map
which sends the class $[\sha]$ of a WKB-algebra to the class $[\sha^+]$ of 
the corresponding WKB-algebroid. 
We refer to \cite{Deligne,BoutetdeMonvel2002} for similar constructions in the 
framework of real manifolds.
\end{remark}

\providecommand{\bysame}{\leavevmode\hbox to3em{\hrulefill}\thinspace}

\end{document}